\documentclass[12pt]{article}

\usepackage{amsmath,amsthm,amsfonts,amssymb}

\newtheorem{theorem}{Theorem}[section]
\newtheorem{corollary}[theorem]{Corollary}
\newtheorem{example}[theorem]{Example}
\newtheorem{proposition}[theorem]{Proposition}

\newtheorem{remark}[theorem]{Remark}

\def\cB{\mathcal{B}}

\def\cF{\mathcal{F}}
\def\cL{\mathcal{L}}

\def\cN{\mathcal{N}}

\def\bE{\mathbb{E}}

\def\bR{\mathbb{R}}

\begin{document}

\title{Intermittency for the wave equation with \\
L\'evy white noise}

\author{Raluca M. Balan\footnote{Corresponding author. Department of Mathematics and Statistics, University of Ottawa,
585 King Edward Avenue, Ottawa, ON, K1N 6N5, Canada. E-mail
address: rbalan@uottawa.ca} \footnote{Research supported by a
grant from the Natural Sciences and Engineering Research Council
of Canada.}\and
Cheikh B. Ndongo\footnote{University of Ottawa, Department of Mathematics and Statistics,
585 King Edward Avenue, Ottawa, ON, K1N 6N5, Canada. E-mail
address: cndon072@uottawa.ca}}

\date{May 15, 2015}
\maketitle

\begin{abstract}
\noindent In this article, we consider the stochastic wave equation on $\bR_{+} \times \bR$ driven by the L\'evy white noise introduced in \cite{B15}. Using Rosenthal's inequality, we develop a maximal inequality for the moments of order $p \geq 2$ of the integral with respect to this noise. Based on this inequality, we show that this equation has a unique solution, which is weakly intermittent in the sense of \cite{FK09,khoshnevisan14}.
\end{abstract}


\vspace{1mm}


\section{Introduction}

In this article, we consider the stochastic wave equation in spatial dimension $d=1$, driven by the L\'evy white noise $L$ introduced in \cite{B15}:
\begin{equation}
\left\{\begin{array}{rcl}
\displaystyle \frac{\partial^2 u}{\partial t^2}(t,x) & = & \displaystyle \frac{\partial^2 u}{\partial x^2}(t,x)+\sigma(u(t,x))\dot{L}(t,x)+b(u(t,x)), \quad t>0, x \in \bR \\[2ex]
\displaystyle u(0,x) & = & v_0(x), \\[1ex]
\displaystyle \frac{\partial u}{\partial t}(0,x) & = & v_1(x),
\end{array}\right. \label{wave} 
\end{equation}
where $v_0$ is a bounded function and $v_1 \in L^1(\bR)$. We assume that $\sigma$ and $b$ are Lipschitz continuous functions. We let $G$ be the fundamental solution of the wave equation on $\bR$:
\begin{equation}
\label{def-G}
G(t,x)=\frac{1}{2}1_{\{|x| \leq t\}}, \quad t>0,x \in \bR,
\end{equation}
and $w$ be the solution of the homogeneous wave equation on $\bR$ with the same initial conditions as \eqref{wave}:
\begin{equation}
\label{def-w}
w(t,x)=\frac{1}{2}\int_{x-t}^{x+t}v_1(y)dy+\frac{1}{2}\Big(v_0(x+t)+
v_0(x-t)\Big).
\end{equation}

We say that a predictable process $u=\{u(t,x); t \geq 0,x \in \bR\}$ is a (mild) {\em solution} of \eqref{wave} if it satisfies the following integral equation:
\begin{eqnarray}
\nonumber
u(t,x)&=&w(t,x)+\int_0^t \int_{\bR}G(t-s,x-y)\sigma(u(s,y))L(ds,dy)+\\
\label{def-sol}
& & \int_0^t \int_{\bR}G(t-s,x-y)b(u(s,y))dyds.
\end{eqnarray}

Before we proceed, we recall briefly from \cite{B15} the definition of the L\'evy white noise $L$ and the construction of the stochastic integral with respect to this noise.

We consider a Poisson random measure (PRM) $N$ on the space $\bE=\bR_{+} \times \bR \times \bR_0$, of intensity $\mu=dtdx \nu(dz)$, where $\bR_0=\bR \verb2\2 \{0\}$ and $\nu$ is a {\em L\'evy measure} on $\bR$, i.e.
$$\int_{\bR_0}(1 \wedge |z|^2)\nu(dz)<\infty \quad \mbox{and} \quad \nu(\{0\})=0.$$
We denote by $\widehat{N}$ the compensated PRM defined by $\widehat{N}(A)=N(A)-\mu(A)$ for any Borel set $A$ in $\bE$ with $\mu(A)<\infty$. Throughout this article, we assume that $\nu$ satisfies:
\begin{equation}
\label{def-v}
m_2:=\int_{\bR_0}|z|^2 \nu(dz)<\infty.
\end{equation}

Suppose that $N$ is defined on a complete probability space $(\Omega, \cF,P)$. On this space, we consider the filtration
$$\cF_t= \sigma(\{N([0,s] \times B \times \Gamma); 0 \leq s \leq t, B \in \cB_b(\bR), \Gamma \in \cB_b(\bR_0)\}) \vee \cN, \quad t \geq 0,$$
where $\cN$ is the class of $P$-negligible sets, $\cB_b(\bR)$ is the class of bounded Borel sets in $\bR$, and $\cB_b(\bR_0)$ is the class of Borel sets in $\bR_0$ which are bounded away from $0$.
Similarly to It\^o's classical theory, for any predictable process $H$ which satisfies
\begin{equation}
\label{cond-H}
E\int_0^t \int_{\bR}\int_{\bR_0}|H(s,x,z)|^2\widehat{N}(ds,dx,dz)<\infty \quad \mbox{for all} \ t>0,
\end{equation}
we can define the stochastic integral of $H$ with respect to $\widehat{N}$, and the integral process \linebreak  $\{\int_0^t \int_{\bR}\int_{\bR_0}H(s,x,z)\widehat{N}(ds,dx,dz); t \geq 0\}$ is a zero-mean square-integrable martingale (see for instance Section 2.2 of \cite{kunita04} or Section 4.2 of \cite{applebaum09}). We work only with c\`adl\`ag modifications of such integral processes. (A process is {\em c\`adl\`ag} if its sample paths are right-continuous and have left limits.) Moreover, the following isometry property holds:
\begin{equation}
\label{E-int-hatN}
E\left|\int_0^t \int_{\bR}\int_{\bR_0}H(s,x,z)\widehat{N}(ds,dx,dz)\right|^2=\int_0^t \int_{\bR}\int_{\bR_0}|H(s,x,z)|^2 \nu(dz)dxds.
\end{equation}
Here, we say that a process $H=\{H(t,x,z); t \geq 0, x \in \bR,z \in \bR_0\}$ is  {\em predictable} if it is measurable with respect to the $\sigma$-field generated by all linear combinations of ``elementary'' processes, i.e. processes of the form
$H(\omega,t,x,z)=Y(\omega)1_{(a,b]}(t)1_{A}(x) 1_{\Gamma}(z)$,
with $0 \leq a<b$, $Y$ an $\cF_a$-measurable bounded random variable, $A \in \cB_b(\bR)$ and $\Gamma \in \cB_b(\bR_0)$.

The {\em L\'evy white noise} defined in \cite{B15} is a ``worthy'' martingale measure $L=\{L_t(B); t \geq 0, B \in \cB_b(\bR)\}$ in the sense of \cite{walsh86}, given by:
$$L_t(B)=\int_0^t \int_{B}\int_{\bR_0}z\widehat{N}(ds,dx,dz).$$
This noise is characterized by the following properties: (i) $L_t(B_1), \ldots, L_t(B_k)$ are independent for any $t>0$ and for any disjoint sets $B_1, \ldots,B_k \in \cB_b(\bR)$; (ii) for any $0 \leq s<t$ and $B \in \cB_b(\bR)$, $L_t(B)-L_s(B)$ is independent of $\cF_s$ and has characteristic function:
$$E(e^{iu(L_t(B)-L_s(B))})=\exp\left\{(t-s)|B|\int_{\bR}(e^{iuz}-1-iuz)\nu(dz)
\right\}, \quad u \in \bR,$$
where $|B|$ is the Lebesgue measure of $B$. Using Walsh' theory developed in \cite{walsh86}, for any predictable process $X=\{X(t,x);t \geq 0,x \in \bR\}$ which satisfies the condition:
\begin{equation}
\label{cond-X}
E\int_0^t \int_{\bR}|X(s,x)|^2 dxds<\infty \quad \mbox{for any} \ t>0,
\end{equation}
we can define the stochastic integral of $X$ with respect to $L$, and the integral process $\{\int_0^t \int_{\bR}X(s,x)L(ds,dx); t \geq 0\}$ is a zero-mean square-integrable martingale. Moreover,
$$\int_0^t \int_{\bR}X(s,x)L(ds,dx)=\int_0^t \int_{\bR}\int_{\bR_0}X(s,x)z\widehat{N}(ds,dx,dz),$$
and therefore, by \eqref{E-int-hatN} and \eqref{def-v},
the following isometry property holds:
\begin{equation}
\label{isometry}
E\left|\int_0^t \int_{\bR}X(s,x)L(ds,dx)\right|^2=m_2 E\int_{0}^{t}\int_{\bR}|X(s,x)|^2dxds.
\end{equation}
This finishes our discussion about the noise $L$.

\vspace{2mm}

We now return to equation \eqref{wave}.
Our goal is to show that, under certain conditions, this equation has a unique solution which is {\em weakly intermittent} in the sense of \cite{FK09, khoshnevisan14}, i.e. for any $x \in \bR$,
$$\underline{\gamma}_x(2)>0 \quad \mbox{and} \quad \overline{\gamma}_x(p)<\infty \quad \mbox{for all} \ p > 2,$$
where $\underline{\gamma}_x(p)$, respectively $\overline{\gamma}_x(p)$ are the {\em lower and upper Lyapunov exponents}, defined by:
$$\underline{\gamma}_x(p)=\liminf_{t \to \infty}\frac{1}{t} \log E|u(t,x)|^p \quad \mbox{and} \quad \overline{\gamma}_x(p)=\limsup_{t \to \infty}\frac{1}{t} \log E|u(t,x)|^p.$$

The following theorem is the main result of this article.

\begin{theorem}
\label{main-th}
Assume that $\sigma$ and $b$ are Lipschitz continuous functions, $v_0$ is bounded, $v_1 \in L^1(\bR)$ and the measure $\nu$ satisfies \eqref{def-v}.

(a) Equation \eqref{wave} has a unique solution $u$ which is continuous in $L^2(\Omega)$ and satisfies:
$$\sup_{(t,x) \in [0,T] \times \bR}E|u(t,x)|^2<\infty \quad \mbox{for all} \ T>0.$$

(b) If the measure $\nu$ satisfies
\begin{equation}
\label{p-moment-nu}
m_p:=\int_{\bR_0}|z|^p \nu(dz)<\infty \quad \mbox{for all} \quad p\geq 2,
 \end{equation}
then for any $t>0$ and $p \geq 2$,
$$\sup_{x \in \bR}E|u(t,x)|^p \leq L_1^p \exp\Big(L_2^{p/2}M_p^{1/2}p^{p/2}t \Big),$$
where $M_p=\max\{m_2,m_p\}$ and $L_1>0,L_2>0$ are constants which do not depend on $p$. Therefore in this case, for any $x \in \bR$ and $p \geq 2$, $\overline{\gamma}_x(p) \leq L_2^{p/2} M_p^{1/2}p^{p/2}$.

(c) If $v_0(x) \geq a>0$ for all $x \in \bR$, $v_1\equiv 0$, $b \equiv 0$ and $L_{\sigma}:=\inf_{x \not=0}|\sigma(x)/x|>0$, then there exists $t_0>0$ such that for any $t\geq t_0$,
$$\inf_{x \in \bR}E|u(t,x)|^2 \geq \frac{a^2}{2} \exp\left\{L_{\sigma} \Big(\frac{m_2}{2}\Big)^{1/2}t\right\}.$$
Therefore in this case, for any $x \in \bR$,
$\underline{\gamma}_x(2) \geq L_{\sigma} (m_2/2)^{1/2}$.
\end{theorem}

\begin{example}[Gamma white noise]
{\rm Recall that the characteristic function of a random variable $X$ with a Gamma distribution with parameters $\alpha>0$ and $\beta>0$ is given by:
$$E(e^{iuX})=\left(\frac{1}{1-iu/\beta} \right)^{\alpha}=\alpha\int_0^{\infty}(e^{iuz}-1) z^{-1}e^{-\beta z}dz.$$
It follows that a L\'evy white noise $L$ with L\'evy measure $\nu(dz)=\alpha z^{-1}e^{-\beta z}1_{\{z>0\}}dz$ can be represented as
$L_t(B)=X_t(B)-E(X_t(B))=X_t(B)-\alpha t|B|/\beta$, where $X_t(B)$ has a Gamma distribution with parameters $\alpha t|B|$ and $\beta$. In this case, we say that $L$ is a {\em Gamma white noise}. Condition  \eqref{p-moment-nu} holds, since $m_p=\alpha \int_0^{\infty}z^{p-1}e^{-\beta z}dz=\alpha \Gamma(p)\beta^{-p}$ for all $p>0$.
}
\end{example}

\begin{remark}
{\rm 
If condition \eqref{p-moment-nu} holds,
it can be proved that
the solution $u$ of \eqref{wave} satisfies:
$E|u(t,x)-u(t',x')|^p \leq C_{p,T,K}(|x-x'|+|t-t'|)$
for any $t,t' \in [0,T]$, $x,x' \in [-K,K]$. 

}
\end{remark}

\begin{remark}
{\rm A zero-mean Gaussian process $W=\{W_t(B);t \geq 0,B \in \cB_b(\bR)\}$ with covariance $E[W_t(A)W_s(B)]=(t \wedge s)|A \cap B|$ is called a {\em Gaussian white noise}. As in \cite{khoshnevisan14}, it can be proved that the solution $u$ of the wave equation with Gaussian white noise is weakly intermittent. In fact,
$\sup_{x \in \bR}E|u(t,x)|^p \leq L_1^p \exp\Big(L_2 p^{3/2}t \Big)$
for any $t>0, p \geq 2$.
}
\end{remark}

\begin{remark}
{\rm We can also consider the heat equation with L\'evy white noise $L$:
\begin{equation}
\label{heat}
\frac{\partial u}{\partial t}(t,x) =  \displaystyle \frac{1}{2} \ \frac{\partial^2 u}{\partial x^2}(t,x)+\sigma(u(t,x))\dot{L}(t,x)+b(u(t,x)), \quad t>0, x \in \bR
\end{equation}
and initial condition $u(0,x)=v_0(x), x \in \bR$, where $\sigma$ and $b$ are Lipschitz functions and $v_0$ is bounded. By classical methods, it follows that equation \eqref{heat}
has a unique solution $u$ which is continuous in $L^2(\Omega)$.
Unfortunately, proving that the solution of this equation 
has finite moments of order $p\geq 3$ (when the L\'evy measure $\nu$ satisfies \eqref{p-moment-nu}) remains an open problem. The difficulty arises due to {\em the second term} in Rosenthal's inequality, which leads to the integral
$\int_{\bR}G^p(s,y)dy=c_p s^{(1-p)/2}=:g(s)$,
where $G(t,x)=(2\pi t)^{-1/2}\exp(-x^2/(2t))$ is the fundamental solution of the heat equation on $\bR$. When $p\geq 3$, $g$ is not integrable around $0$ and the standard argument based on Gronwall lemma cannot be applied.
}
\end{remark}

This article is organized as follows. In Section \ref{max-ineq-section}, we prove a Rosenthal-type inequality for the moments of order $p \geq 2$ of the stochastic integral with respect to the noise $L$. Section \ref{proof-section} is dedicated to the proof of Theorem \ref{main-th}. Classical methods yield immediately the existence and uniqueness of the solution. The major effort will be dedicated to proving the upper bound. For this, we follow the approach presented in Chapter 5 of \cite{khoshnevisan14} for the heat equation with Gaussian white noise, replacing the
Burkholder-Davis-Gundy inequality by the Rosenthal-type inequality mentioned above. Luckily, due to the isometry property \eqref{isometry}, the lower bound requires very little effort, using the same argument as in
Chapter 7 of \cite{khoshnevisan14} in the case of the heat equation with Gaussian white noise.

\section{Maximal Inequality}
\label{max-ineq-section}

In this section, we give a maximal inequality for the moments of order $p \geq 2$ of the stochastic integrals with respect to $\widehat{N}$ and $L$. We denote by $\|\cdot\|_p$ the norm in $L^p(\Omega)$.

We begin by recalling Rosenthal's inequality for c\`adl\`ag martingales (see e.g. Lemma 2.1 of \cite{DV90}).

\begin{theorem}[Rosenthal's inequality]
\label{DV-theorem}
Let $\{M(t)\}_{t \geq 0}$ be a c\`adl\`ag square-integrable martingale with $M(0)=0$ and $\langle M \rangle$ be its predictable quadratic variation. We denote by $(\Delta M)(t)$ the jump size of $M$ at time $t$. Then for any $p \geq 2$, there exists a constant $B_p>0$ depending on $p$ such that for any $t>0$,
$$\|\sup_{s \leq t}|M(s)|\|_p \leq B_p \Big(\| \langle M \rangle(t)^{1/2}\|_p +\|\sup_{s \leq t}|(\Delta M)(s)| \|_p\Big).$$
\end{theorem}

\begin{remark}
{\rm By Theorem 1 of \cite{ren-tian03}, we know that the constant $B_p$ in Rosenthal's inequality has the following asymptotic behaviour:
\begin{equation}
\label{constant-Bp}
B_p=O\left(\frac{p}{\ln p}\right), \quad \mbox{as} \ p \to \infty.
\end{equation}
 }
\end{remark}

Based on Rosenthal's inequality, we obtain the following maximal inequality for the stochastic integral with respect to $\widehat{N}$.

\begin{theorem}
\label{Rosenthal-ineq-th}
Let $Y=\{Y(t)\}_{t \geq 0}$ be a process given by:
$$Y(t)=\int_0^t \int_{\bR}\int_{\bR_0}H(s,x,z)\widehat{N}(ds,dx,dz),$$
where $H$ is a predictable process which satisfies \eqref{cond-H}. Then, there exists a c\`adl\`ag modification of $Y$ (denoted also by $Y$) such that for any $t>0$ and for any $p \geq 2$,
\begin{eqnarray}
\nonumber
\left\|\sup_{s \leq t}|Y(s)|\right\|_p
& \leq & B_p \left\{ \left\| \left(\int_0^t \int_{\bR}\int_{\bR_0}
|H(s,x,z)|^2 \nu(dz)dxds \right)^{1/2} \right\|_{p} \right. \\
\label{Rosenthal-ineq-N}
& & + \left. \left(E\int_0^t \int_{\bR}\int_{\bR_0}|H(s,x,z)|^p \nu(dz)dxds \right)^{1/p}\right\},
\end{eqnarray}
where $B_p$ is the constant in Rosenthal's inequality.
\end{theorem}

\noindent {\bf Proof:}
{\em Case 1.} We assume that the process $Y$ is of the form:
\begin{equation}
\label{case1}
Y(t)=\int_0^t \int_{B}\int_{\{|z|>\varepsilon\}}H(s,x,z)\widehat{N}(ds,dx,dz)
\end{equation}
for some $B \in \cB_b(\bR)$ and $\varepsilon>0$. In this case, $Y(t)=Y_d(t)+Y_c(t)$, where
\begin{eqnarray*}
Y_d(t) &=& \int_0^t \int_{B}\int_{\{|z|>\varepsilon\}}H(s,x,z)N(ds,dx,dz),\\
Y_c(t) &=&
\int_0^t \int_{B}\int_{\{|z|>\varepsilon\}}H(s,x,z)\nu(dz)dxds
\end{eqnarray*}
Note that $Y_d$ is c\`adl\`ag and $Y_c$ is continuous. Hence, $Y$ is a c\`adl\`ag martingale. Its predictable quadratic variation is given by:
$$\langle Y \rangle (t)=\int_0^t \int_{B}\int_{\{|z|>\varepsilon\}}|H(s,x,z)|^2 \nu(dz)dxds.$$

Let $\Gamma=\{z;|z|>\varepsilon\}$ and $\lambda=|B|\nu(\Gamma)$. By Proposition 5.3 of \cite{resnick07}, we may assume that the points of $N$ in $\bR_{+} \times B \times \Gamma$ are given by $\{(T_i,X_i,Z_i)\}_{i \geq 1}$, where $T_1<T_2< \ldots$ are the points of a Poisson process on $\bR_{+}$ of intensity $\lambda$, and $\{(X_i,Z_i)\}_{i \geq 1}$ are i.i.d. random variables on $B \times \Gamma$ with law $dx \nu(dz)/\lambda$. Hence, the sample path $t \mapsto Y_d(t)$ is a step function given by:
$$Y_d(t)=\sum_{i=1}^{n}H(T_i,X_i,Z_i) \quad \mbox{if} \quad T_{n-1} \leq t<T_n.$$
We infer that $Y_d$ has a jump of size $H(T_i,X_i,Z_i)$ at $T_i$, and there are a finite number of such jumps in the interval $[0,t]$.
The conclusion follows by applying Theorem \ref{DV-theorem} to the c\`adl\`ag martingale $M=Y$, since
\begin{eqnarray*}
\lefteqn{E \left(\sup_{s \leq t}|\Delta Y(s)|^p \right) = E \left(\sup_{s \leq t}|\Delta Y_d(s)|^p \right) \leq E \left(\sum_{s \leq t}|\Delta Y_d(s)|^p \right)=} \\
& &  E \left(\sum_{T_i \leq t}|H(T_i,X_i,Z_i)|^p \right)
=E \int_0^t \int_B \int_{\{|z|>\varepsilon\}}|H(s,x,z)|^p N(ds,dx,dz)=\\
& & E \int_0^t \int_B \int_{\{|z|>\varepsilon\}}|H(s,x,z)|^p \nu(dz)dxds.
\end{eqnarray*}

{\em Case 2.} We assume that $H(s,x,z)=0$ if $|z| \leq \varepsilon$, for some $\varepsilon>0$. We suppose that the right-hand side of \eqref{Rosenthal-ineq-N} is finite; otherwise, the result is trivial.
Let $Y_k$ be the c\`adl\`ag process given by \eqref{case1} with $B=E_k$, where $(E_k)_{k \geq 1}$ is a sequence in $\cB_b(\bR)$ such that $E_k \subset E_{k+1}$ for all $k$ and $\bigcup_{k \geq 1} E_k=\bR$.
We apply the result of {\em Case 1} to the process $Y_k-Y_l$, for $k>l$. Hence,
$$\left\|\sup_{s \leq t}|Y_k(s)-Y_l(s)|  \right\|_{p} \leq
B_p (A_{k,l}+B_{k,l}),$$
where
\begin{eqnarray*}
A_{k,l}&=& \left\| \left(\int_0^t \int_{E_k- E_l}\int_{\{|z|>\varepsilon\}}
|H(s,x,z)|^2 \nu(dz)dxds \right)^{1/2} \right\|_{p}  \\
B_{k,l} & =&  \left(E\int_0^t \int_{E_k -E_l}\int_{\{|z|>\varepsilon\}}|H(s,x,z)|^p \nu(dz)dxds \right)^{1/p}.
\end{eqnarray*}
By the dominated convergence theorem, $A_{k,l}\to 0$ and $B_{k,l} \to 0$ as $k,l \to \infty$. Hence, $(Y_k)_{k \geq 1}$ is a Cauchy sequence in $L^p(\Omega; D[0,t])$, where $D[0,t]$ is the set of c\`adl\`ag functions on $[0,t]$ equipped with the sup norm. Its limit is a c\`adl\`ag modification of $Y$. The conclusion follows by applying the result of {\em Case 1} to the process $Y_k$, and letting $k \to \infty$.

{\em Case 3.} In the general case, we consider a sequence $(\varepsilon_k)_{k \geq 1}$ with $\varepsilon_k \downarrow 0$ and we let $Y_k$ be the process given by \eqref{case1} with $B=\bR$ and $\varepsilon=\varepsilon_k$. We denote also by $Y_k$ the c\`adl\`ag modification of this process given by {\em Case 2}. We apply the result of {\em Case 2} to $Y_k-Y_l$ for $k>l$. As above, we infer that $(Y_k)_{k \geq 1}$ is a Cauchy sequence in $L^p(\Omega; D[0,t])$. Its limit is a c\`adl\`ag modification of $Y$ for which the conclusion holds. $\Box$

\vspace{3mm}

The following result plays an important role in the proof of Theorem \ref{main-th}.

\begin{corollary}
\label{Rosenthal-ineq-cor}
Let $Y=\{Y(t)\}_{t \geq 0}$ be a process given by $$Y(t)=\int_0^t\int_{\bR}X(s,x)L(ds,dx),$$ where $X=\{X(t,x);t \geq 0,x \in \bR\}$ is a predictable process which satisfies \eqref{cond-X}.
Suppose that the measure $\nu$ satisfies condition \eqref{p-moment-nu}.
Then there exists a c\`adl\`ag modification of $Y$ (denoted also by $Y$) such that for any $t>0$ and for any $p\geq 2$,
\begin{eqnarray}
\nonumber
\left\|\sup_{s \leq t}|Y(s)|\right\|_{p} & \leq & B_p \left\{m_2^{1/2}\left\| \left(\int_0^t \int_{\bR} |X(s,x)|^2 dxds \right)^{1/2}\right\|_{p}+ \right. \\
\label{Rosenthal-ineq-cor1}
& & \quad \quad \left. m_p^{1/p}\left(E\int_0^t \int_{\bR} |X(s,x)|^p dxds\right)^{1/p} \right\},
\end{eqnarray}
where $B_p$ is the constant in Rosenthal's inequality.
\end{corollary}

\noindent {\bf Proof:} The result follows by applying Theorem \ref{Rosenthal-ineq-th} with $H(s,x,z)=X(s,x)z$. $\Box$

\section{Proof of Theorem \ref{main-th}}
\label{proof-section}

This section is dedicated to the proof of Theorem \ref{main-th}.

\vspace{3mm}

As in \cite{FK09, khoshnevisan14}, for any $\beta>0$ and $p \geq 2$, we consider the space $\cL^{\beta,p}$ of predictable processes $\Phi=\{\Phi(t,x); t \geq 0,x \in \bR\}$ such that
$\|\Phi\|_{\beta,p}<\infty$, where
$$\|\Phi\|_{\beta,p}:=\sup_{t \geq 0} \sup_{x \in \bR}\
\Big(e^{-\beta t}\|\Phi(t,x)\|_p\Big).$$
It can be proved that $\cL^{\beta,p}$ equipped with $\|\cdot\|_{\beta,p}$ is a Banach space. (We identify processes $\Phi_1$ and $\Phi_2$, if $\Phi_1$ is a modification of $\Phi_2$.)

We begin with some preliminary results.

The first result is a variant of the Rosenthal's inequality, which gives the analogue of Proposition 4.4 of \cite{khoshnevisan14} for the L\'evy white noise $L$.

\begin{proposition}
\label{Rosenthal-ineq-prop}
For any measurable function $h: \bR_{+} \times \bR \to \bR$ such that $h \in L^2([0,t] \times \bR)$ for all $t>0$ and for any $\Phi \in \cup_{\beta>0}\cL^{\beta,2}$, the process $\{\int_0^t \int_{\bR} h(s,y)\Phi(s,y)L(ds,dy), t\geq 0\}$ is well-defined and has a c\`adl\`ag modification which satisfies:
\begin{eqnarray*}
\left\|\int_0^t \int_{\bR}h(s,y) \Phi(s,y)L(ds,dy)\right\|_p  & \leq & B_p \left\{m_2^{1/2} \left(\int_0^t \int_{\bR}|h(s,y)|^2 \|\Phi(s,y)\|_{p}^2 dyds \right)^{1/2} \right. \\
& &  + \left. m_p^{1/p} \left(\int_0^t \int_{\bR}|h(s,y)|^p \|\Phi(s,y)\|_p^p dyds \right)^{1/p} \right\}.
\end{eqnarray*}
\end{proposition}

\noindent {\bf Proof:} To see that the stochastic integral is well-defined, we note that:
$$E\int_0^t \int_{\bR}|h(s,y)|^2 |\Phi(s,y)|^2 dyds \leq \|\Phi\|_{\beta,2}^2 e^{2\beta t}\int_{0}^{t} \int_{\bR}|h(s,y)|^2dyds<\infty.$$
We apply Corollary \ref{Rosenthal-ineq-cor} with $X(s,y)=h(s,y)\Phi(s,y)$. Then, we use Minkowski's inequality with respect to the norm $\|\cdot\|_{p/2}$ for the first term on the right-hand side of \eqref{Rosenthal-ineq-cor1}. $\Box$

\vspace{3mm}

For any $\Phi \in \cup_{\beta>0}\cL^{\beta,2}$ we define the {\em stochastic convolution} with $G$ by:
$$(G * \Phi)(t,x)=\int_0^t \int_{\bR}G(t-s,x-y)\Phi(s,y)dyds$$
for all $t>0$ and $x \in \bR$, and $(G * \Phi)(0,x)=0$ for all $x \in \bR$.
Since $G(t-\cdot,x-\cdot) \in L^2([0,t] \times \bR$ for any $t>0$ and $x \in \bR$, $G * \Phi$ is well-defined by Proposition \ref{Rosenthal-ineq-prop}.

\vspace{1mm}

The following result is the analogue of Proposition 5.2 of \cite{khoshnevisan14} for $L$.

\begin{proposition}[Stochastic Young Inequality for the L\'evy white noise]
\label{Young-prop}
For any $\beta>0$, $p \geq 2$ and $\Phi \in \cL^{\beta,2}$,
$$\|G * \Phi\|_{\beta,p} \leq B_p \left\{\frac{m_2^{1/2}}{2\sqrt{2} \beta}+\frac{1}{2}\left(\frac{2m_p}{p^2 \beta^2} \right)^{1/p} \right\}\|\Phi\|_{\beta,p}.$$
\end{proposition}

\noindent {\bf Proof:} Fix $(t,x)$. We apply Proposition \ref{Rosenthal-ineq-prop} with $h(s,y)=G(t-s,x-y)$. We assume that $\|\Phi\|_{\beta,p}<\infty$; otherwise the result is trivial.
Since $\|\Phi(s,y)\|_p \leq e^{\beta s}\|\Phi\|_{\beta,p}$, we obtain:
\begin{eqnarray*}
\|(G* \Phi)(t,x)\|_p & \leq & B_p \|\Phi\|_{\beta,p} \left\{m_2^{1/2} \left( \int_0^t e^{2\beta s} \int_{\bR}G^2(t-s,x-y)dyds\right)^{1/2}\right.\\
& & \left. + m_p^{1/p} \left(\int_0^t e^{p\beta s}\int_{\bR}G^p(t-s,x-y)dyds \right)^{1/p} \right\}.
\end{eqnarray*}

\noindent By the definition \eqref{def-G} of $G$, we have $\int_{\bR}G^p(t-s,x-y)dy=2^{1-p}s$. Hence, for any $p \geq 2$,
\begin{eqnarray*}
& & \int_0^t e^{p\beta s} \int_{\bR}G^p(t-s,x-y)dyds=2^{1-p}\int_0^t (t-s)e^{p\beta s}ds\leq \\
& & \quad \quad 2^{1-p}e^{p\beta t}\int_0^{\infty} se^{-p\beta s}ds=2^{1-p}e^{p\beta t}\frac{1}{(p\beta)^2}.
\end{eqnarray*}
It follows that
$$\|(G* \Phi)(t,x)\|_p  \leq B_p \|\Phi\|_{\beta,p} e^{\beta t}
\left\{\frac{m_2
^{1/2}}{2\sqrt{2}\beta} +\frac{1}{2} \left(\frac{2m_p}{p^2 \beta^2} \right)^{1/p} \right\}.$$
We multiply by $e^{-\beta t}$ and we take the supremum over $t \geq 0$ and $x \in \bR$. $\Box$

\begin{remark}
{\rm In the case of the stochastic convolution with the space-time Gaussian white noise $W$, we have the following inequality:
$$\left\| \int_0^t \int_{\bR}G(t-s,x-y)\Phi(s,y)W(ds,dy)\right\|_{\beta,p} \leq z_p \frac{1}{2\sqrt{2}\beta}\|\Phi\|_{\beta,p},$$
where $z_p \sim 2\sqrt{p} $ is the constant in Burkholder-Davis-Gundy inequality.
}
\end{remark}

We now return to the study of equation \eqref{wave}. Since the functions $\sigma$ and $b$ are Lipschitz continuous, there exists a constant $L>0$ such that
\begin{equation}
\label{Lip1}
\max\{|\sigma(x)-\sigma(y)|,|b(x)-b(y)|\} \leq L|x-y| \quad \mbox{for all} \ x,y \in \bR.
\end{equation}
Choosing $L>\max\{|\sigma(0)|,|b(0)|\}$, we have
\begin{equation}
\label{Lip2}
\max\{|\sigma(x)|,|b(x)|\} \leq L(1+|x|) \quad \mbox{for all} \ x \in \bR.
\end{equation}

We consider the sequence $(u_n)_{n \geq 0}$ of Picard iterations, defined by:
\begin{eqnarray}
\nonumber
u_{n+1}(t,x)&=& w(t,x)+\int_0^t \int_{\bR}G(t-s,x-y)\sigma(u_n(s,y))L(ds,dy)+\\
\label{def-Picard}
& &
\int_0^t \int_{\bR}G(t-s,x-y)b(u_n(s,y))duds,
\end{eqnarray}
for all $n \geq 0$, and $u_0(t,x)=w(t,x)$.

Similarly to Proposition 5.8 of \cite{khoshnevisan14}, we prove the following result.

\begin{proposition}
The processes $(u_n)_{n \geq 0}$ are well-defined and belong to $\cup_{\beta>0}\cL^{\beta,2}$.
If condition \eqref{p-moment-nu} holds, then there exist some constants $L_1>0$ and $L_2>0$ such that
\begin{equation}
\label{Picard-ineq}
\sup_{x \in \bR}E|u_n(t,x)|^p \leq L_1^p \exp\Big(L_2^{p/2} p^{p/2}M_p^{1/2}t\Big),
\end{equation}
for any $t>0,p \geq 2$ and $n \geq 0$, where $M_p=\max\{m_2,m_p\}$.
\end{proposition}

\noindent {\bf Proof:} The fact that the processes $(u_n)_{n \geq 1}$ are well-defined follows by classical methods. More precisely, by induction on $n$, it can be proved that: {\em (i)} $u_n(t,x)$ is well-defined; {\em (ii)} $\sup_{(t,x) \in [0,T] \times \bR}E|u_n(t,x)|^2<\infty$ for any $T>0$; {\em (iii)} the map $(t,x) \mapsto u_n(t,x)$ is continuous in $L^2(\Omega)$; {\em (iv)} $u_n \in \cL^{\beta,2}$ for all $\beta>0$. At each step, we work with a predictable modification of $u_n$ (denoted also by $u_n$). We omit the details.

It remains to prove \eqref{Picard-ineq}. We will estimate $E|u_{n+1}(t,x)|^p$ by treating separately the three terms on the right-hand side of \eqref{def-Picard}. For the first term, recalling the definition \eqref{def-w} of $w$, we see that 
$|w(t,x)|\leq \frac{1}{2}\int_{\bR}|v_1(x)|dx+\sup_{x \in \bR}|v_0(x)|=:K$ for all $t>0$ and $x \in  \bR$, and hence
\begin{equation}
\label{bound-w}
\|w\|_{\beta,p} \leq K.
\end{equation}

For the second term, by Proposition \ref{Young-prop} and property \eqref{Lip2}, we have:
\begin{equation}
\label{bound-Bn}
\|G * \sigma(u_n)\|_{\beta,p} \leq C_{\beta,p}\|\sigma(u_n)\|_{\beta,p}
 \leq   C_{\beta,p}
L \Big(1+\|u_n\|_{\beta,p}\Big),
\end{equation}
where
\begin{equation}
\label{def-Cbp}
C_{\beta,p}=B_p\left\{\frac{m_2^{1/2}}{2\sqrt{2} \beta}+\frac{1}{2}\left(\frac{2m_p}{p^2 \beta^2} \right)^{1/p}\right\}.
\end{equation}

We denote by $C_n(t,x)$ the third term on the right-hand side of \eqref{def-Picard}. By Minkowski's inequality, property \eqref{Lip2}, and the fact that $\int_{\bR}G(t-s,x-y)dy=t-s$, we have:
\begin{eqnarray*}
\lefteqn{\|C_n(t,x)\|_p \leq  L \int_0^t \int_{\bR}G(t-s,x-y)(1+\|u_n(s,y)\|_p) dyds \leq } \\
& &  L \int_0^t \Big(1+\sup_{y \in \bR}\|u_n(s,y)\|_p\Big)(t-s)ds \leq L \left\{\frac{t^2}{2}+\|u_n\|_{\beta,p}\int_0^t e^{\beta s}(t-s)ds \right\} \leq \\
& &  L \left\{ \frac{t^2}{2}+\|u_n\|_{\beta,p} e^{\beta t} \int_0^{\infty} s e^{-\beta s}ds\right\}= L \left\{ \frac{t^2}{2}+\|u_n\|_{\beta,p} e^{\beta t} \frac{1}{\beta^2}\right\}.
\end{eqnarray*}

\noindent We multiply by $e^{-\beta t}$ and we take the supremum over $t \geq 0$ and $x \in \bR$. Using the fact that $\sup_{t \geq 0}t^2 e^{-\beta t}=4/(e \beta)^2$, we obtain:
\begin{equation}
\label{bound-Cn}
\|C_n\|_{\beta,p} \leq \frac{2L}{e^2 \beta^2}+\frac{L}{\beta^2}\|u_n\|_{\beta,p}.
\end{equation}

Using Minkowski's inequality with respect to the norm $\|\cdot\|_{\beta,p}$ and relations \eqref{bound-w}, \eqref{bound-Bn} and \eqref{bound-Cn}, we obtain: \begin{equation}
\label{bound-un}
\|u_{n+1}\|_{\beta,p} \leq K+C_{\beta,p}L+\frac{2L}{e^2\beta^2}+L\left(C_{\beta,p}+\frac{1}{\beta^2}\right)
\|u_n\|_{\beta,p}.
\end{equation}
By induction on $n$, it follows that $\|u_n\|_{\beta,p}<\infty$ for all $\beta>0$ and $p \geq 2$.

We now prove \eqref{Picard-ineq}. We fix $p \geq 2$. We would like to choose $\beta>0$ (depending on $p$) such that
\begin{equation}
\label{cond-beta}
\frac{L}{\beta^2}<\frac{1}{4} \quad \mbox{and} \quad C_{\beta,p}L<\frac{1}{4}.
\end{equation}
Using definition \eqref{def-Cbp} of $C_{\beta,p}$, we see that it suffices to pick $\beta>0$ such that
$$\frac{L}{\beta^2}<\frac{1}{4}, \quad B_p L \frac{m_2^{1/2}}{2\sqrt{2}\beta}<\frac{1}{8} \quad \mbox{and} \quad \frac{1}{2}B_p L \left(\frac{2m_p}{p^2 \beta^2}\right)^{1/p}<\frac{1}{8}.$$

\noindent By \eqref{constant-Bp}, there exists a constant $C_0 \geq 1$ such that $B_p \leq C_0p$ for all $p \geq 2$. Hence, it is enough to choose $\beta>0$ such that
$$\beta^2>4L, \quad \beta^2>\frac{m_2}{2} (4LC_0 p)^2 \quad \mbox{and} \quad \beta^2>\frac{2m_p}{p^2}(4LC_0 p)^p.$$
We may assume that $L>1/(4m_2)$, and so $m_2 (4LC_0 p)^2/2 \geq m_2 (4L)^2 \geq 4L$. Note that $2p^{p-2} \geq p^2/2$
and $(4LC_0)^p \geq (4LC_0)^2$ (if $4L\geq 1$). Hence, it is enough to choose $\beta>0$ such that $\beta^2>2  p^{p-2} M_p (4C_0L)^p$, where $M_p=\max\{m_2,m_p\}$. We let $\beta>0$ be given by:
\begin{equation}
\label{def-beta}
\beta^2=p^{p-2}M_p L_2^p \quad \mbox{where} \quad L_2=9C_0L.
\end{equation}
With this choice of $\beta$, condition \eqref{cond-beta} is satisfied.

From \eqref{bound-un} and \eqref{cond-beta}, we obtain that
$\|u_{n+1}\|_{\beta,p} \leq \gamma+\frac{1}{2}\|u_n\|_{\beta,p}$,
where $\gamma=K+1/4+1/(2e^2)$. Using recursively this inequality, it follows that
\begin{equation}
\label{bound-un1}
\|u_{n+1}\|_{\beta,p} \leq \sum_{j=0}^{n}\frac{\gamma}{2^j}+\frac{1}{2^{n+1}}\|u_0\|_{\beta,p} \leq 2\gamma+\frac{1}{2}K=:L_1 \quad \mbox{for all} \ n \geq 0.
\end{equation}
Hence $E|u_{n+1}(t,x)|^p \leq L_1^p e^{p \beta t}$ and the conclusion follows by the choice of $\beta$. $\Box$

\vspace{3mm}

\noindent {\bf Proof of Theorem \ref{main-th}:}
(a) By classical methods, it can be proved that the sequence $(u_n)_{n \geq 0}$ converges in $L^2(\Omega)$ uniformly in $(t,x) \in [0,T] \times \bR$, and its limit is the solution of equation \eqref{wave} (see for instance the proof of Theorem 13 of \cite{dalang99}). We omit the details.

(b) We proceed as in the proof of Theorem 5.5 of \cite{khoshnevisan14}. We have
\begin{eqnarray*}
\lefteqn{
u_{n+1}(t,x)-u_n(t,x)=\int_0^t \int_{\bR}G(t-s,x-y)[\sigma(u_n(s,y))-\sigma(u_{n-1}(s,y))]L(ds,dy)+}\\
& & \int_0^t \int_{\bR}G(t-s,x-y)[b(u_n(s,y))-b(u_{n-1}(s,y))]dyds=:J_1(t,x)+J_2(t,x)
\end{eqnarray*}

By Proposition \ref{Young-prop} and property \eqref{Lip1}, we have:
$$\|J_1\|_{\beta,p}\leq C_{\beta,p}\|\sigma(u_n)-\sigma(u_{n-1})\|_{\beta,p} \leq C_{\beta,p}L\|u_n-u_{n-1}\|_{\beta,p},$$
where $C_{\beta,p}$ is given by \eqref{def-Cbp}.
Similarly to \eqref{bound-Cn}, it can be shown that
$$\|J_2\|_{\beta,p} \leq \frac{L}{\beta^2}\|u_n-u_{n-1}\|_{\beta,p}.$$
It follows that for any $p \geq 2$, $\beta>0$ and $n \geq 0$,
$$\|u_{n+1}-u_n\|_{\beta,p} \leq L \Big(C_{\beta,p}+\frac{1}{\beta^2}\Big)\|u_n-u_{n-1}\|_{\beta,p}.$$

We fix $p \geq 2$ and let $\beta$ be given by \eqref{def-beta}. (Note that {\em $\beta$ depends on $p$.}) By \eqref{cond-beta},
$$\|u_{n+1}-u_n\|_{\beta,p} \leq \frac{1}{2}\|u_{n}-u_{n-1}\|_{\beta,p} \quad \mbox{for all} \quad n \geq 0.$$
Using this inequality recursively and \eqref{bound-un1}, we obtain that
$$\|u_{n+1}-u_n\|_{\beta,p} \leq \frac{1}{2^n}\|u_1-u_0\|_{\beta,p}\leq \frac{1}{2^{n}}(2L_1) \quad \mbox{for all} \ n \geq 0.$$
It follows that $(u_n)_{n \geq 0}$ is a Cauchy sequence in $\cL^{\beta,p}$. Its limit in this space is the solution $u$ of equation \eqref{wave}. By \eqref{bound-un1}, $\|u\|_{\beta,p} \leq L_1$. Hence, using definition \eqref{def-beta} of $\beta$, we have:
$$E|u(t,x)|^p \leq L_1^p e^{p \beta t}=L_1^p \exp\Big(p^{p/2}M_p^{1/2}L_2^{p/2}t \Big).$$

(c) We argue as in the proof of Theorem 7.8 of \cite{khoshnevisan14}. By \eqref{def-sol}, $E(u(t,x))=w(t,x)$. (Recall that the stochastic integral with respect to $L$ has mean zero and $b \equiv 0$.)
Combined with the isometry property \eqref{isometry}, this yields:
$$E|u(t,x)|^2-w(t,x)^2=E|u(t,x)-w(t,x)|^2=
m_2 E\int_0^t \int_{\bR}G^2(t-s,x-y)|\sigma(u(s,y))|^2dyds.$$
By our hypotheses on $v_0$ and $v_1$, $w(t,x) \geq a$ for all $t>0$ and $x \in \bR$. Since $|\sigma(x)|\geq L_{\sigma}|x|$ for all $x \in \bR$, we infer that:
$$E|u(t,x)|^2 \geq a^2+ L_{\sigma}^2 m_2 \int_0^t \int_{\bR}G^2(t-s,x-y)E|u(s,y)|^2dyds$$
for any $t>0$ and $x \in \bR$. We denote $I(t)=\inf_{x \in \bR}E|u(t,x)|^2$. Using the fact that $\int_{\bR}G^2(t-s,x-y)dy=(t-s)/2$, we obtain that:
$$I(t) \geq a^2+\frac{L_{\sigma}^2 m_2}{2}\int_0^t I(s)(t-s)ds.$$
By Theorem 7.11 of \cite{khoshnevisan14}, it follows that
$I(t) \geq f(t)$ for all $t>0$, where $f$ is a solution of
$$f(t)=a^2+\int_0^t f(s)g(t-s)ds$$ with
$g(t)=(L_{\sigma}^2 m_2/2) t$. Multiplying by $e^{-\lambda t}$ with $\lambda>0$, we get:
$$\widetilde{f}(t)=a^2 e^{-\lambda t}+\int_0^t \widetilde{f}(t-s)\widetilde{g}(s)ds,$$
where $\widetilde{h}(t)=e^{-\lambda t}h(t)$ for any function $h$.
Choosing $\lambda^2=L_{\sigma}^2 m_2/2$, $\widetilde{g}$ becomes a density function on $(0,\infty)$. With this choice of $\lambda$, using Feller's renewal theory, we infer that
$$\lim_{t \to \infty}\widetilde{f}(t)=\frac{\int_0^{\infty}a^2 e^{-\lambda t}dt}{\int_0^{\infty}t\widetilde{g}(t)dt}=\frac{a^2\lambda^{-1}}{m_2 L_{\sigma}^2 \lambda^{-3}}=\frac{a^2}{m_2 L_{\sigma}^2}\lambda^2=\frac{a^2}{2}.$$
Hence,
$\liminf_{t \to \infty}(e^{-\lambda t}I(t)) \geq a^2/2$, i.e. $I(t) \geq (a^2/2)e^{\lambda t}$ for all $t \geq t_0$. $\Box$

\end{document}